\documentclass[11pt]{article}

\usepackage{amsfonts,amssymb,amsmath,latexsym,color,epsfig}
\newtheorem{theorem}{Theorem}
\newtheorem{lemma}[theorem]{Lemma}

\newtheorem{conjecture}{Conjecture}
\newtheorem{corollary}[theorem]{Corollary}

\def\qed{\ifvmode\mbox{ }\else\unskip\fi\hskip 1em plus 10fill$\Box$}
\input{epsf}

\makeatletter
\def\Ddots{\mathinner{\mkern1mu\raise\p@
\vbox{\kern7\p@\hbox{.}}\mkern2mu
\raise4\p@\hbox{.}\mkern2mu\raise7\p@\hbox{.}\mkern1mu}}
\makeatother

\title{\vspace{-0.7cm}Applications of a new separator theorem for string graphs}
\author{
Jacob Fox\thanks{
    Department of Mathematics,
    Massachusetts Institute of Technology,
    Cambridge, MA 02139-4307.
    Email: {\tt fox@math.mit.edu}.
    Research supported by a Simons Fellowship, by NSF grant DMS-1069197, by an Alfred P. Sloan Fellowship, and by an MIT NEC Corporation Award.}  \and J\'anos Pach\thanks{EPFL, Lausanne and Courant Institute, New York, NY. Supported by Hungarian Science Foundation EuroGIGA Grant OTKA NN 102029, by Swiss National Science Foundation Grants 200020-144531 and 200021-137574, and by NSF Grant CCF-08-30272. Email: {\tt  pach@cims.nyu.edu}.}
}
\date{}
\begin{document}
\maketitle

\begin{abstract}
An intersection graph of curves in the plane is called a {\em string graph}. Matou\v sek almost completely settled a conjecture of the authors by showing that every string graph of $m$ edges admits a vertex separator of size $O(\sqrt{m}\log m)$. In the present note, this bound is combined with a result of the authors, according to which every dense string graph contains a large complete balanced bipartite graph. Three applications are given concerning string graphs $G$ with $n$ vertices: (i) if $K_t\not\subseteq G$ for some $t$, then the chromatic number of $G$ is at most $(\log n)^{O(\log t)}$; (ii) if $K_{t,t}\not\subseteq G$, then $G$ has at most $t(\log t)^{O(1)}n$ edges,; and (iii) a lopsided Ramsey-type result, which shows that the Erd\H os-Hajnal conjecture almost holds for string graphs.
\end{abstract}

\section{Introduction}

A graph $G=(V,E)$ is called a {\it string graph} if it is the intersection graph of curves in the plane, i.e., if there is a collection of curves (``strings") $\gamma_v$ in the plane, one curve for each vertex $v \in V$, such that two curves $\gamma_u$ and $\gamma_v$ intersect if and only if $u$ and $v$ are adjacent in $G$.

A {\em separator} in a graph $G=(V,E)$ is a subset $S$ of the vertex set $V$ such that no connected component of $G\setminus S$  has more than $\frac{2}{3}|V|$ vertices. Equivalently, $S$ is a separator of $G$ if there is a partition $V=S \cup V_1 \cup V_2$ with $|V_1|,|V_2| \leq \frac{2}{3}|V|$ such that no vertex in $V_1$ is adjacent to any vertex in $V_2$.

In \cite{FP10}, we proved that every string graph $G$ with $m$ edges has a separator of size $O(m^{3/4}\sqrt{\log m})$, and conjectured that this bound can be improved to $O(\sqrt{m})$.\footnote{Throughout this paper, all logarithms are base 2. Also, for the sake of simplicity, we systematically omit floor and ceiling signs whenever they are not crucial.} This result, if true, would be be best possible. In~\cite{FP08}, we proved our conjecture in the special case where the vertices of $G$ can be represented by curves in the plane with the property that every pair of them intersect in at most a bounded number of points. The starting point of our investigations was a recent paper of Matou\v sek \cite{M13}, in which he ingeniously adapted some powerful techniques developed by Feige, Hajiaghayi, and Lee~\cite{FHL}, using the framework of multicommodity flows to design efficient approximation algorithms for finding small separators in general graphs. (See \cite{BLR}, for a similar application.) Matou\v sek~\cite{M13} proved our above conjecture up to a logarithmic factor.

\begin{lemma}\label{lemsep} {\rm \cite{M13}}
Every string graph with $m$ edges has a separator of size at most $d\sqrt{m}\log m$, where $d$ is an absolute constant.
\end{lemma}

The aim of this note is to combine Lemma~\ref{lemsep} with some previous results of the authors to substantially improve the best known estimates for various important parameters of string graphs.

Our first result provides an upper bound on the chromatic number of string graphs with no complete subgraph of size $t$, which is polylogarithmic in the number of vertices.

\begin{theorem}\label{thm1}
There is an absolute constant $C$ such that every $K_t$-free string graph on $n$ vertices has chromatic number at most $(\log n)^{C \log t}$.
\end{theorem}

Previously, it was not even known if the chromatic number of every triangle-free string graph on $n$ vertices is at most $n^{o(1)}$. In the other direction, solving an old problem of Erd\H os, for every $n$, Pawlik, Kozik, Krawczyk, Laso\'n, Micek, Trotter, and Walczak~\cite{PKKLMTW} constructed a triangle-free intersection graph of $n$ {\em segments} in the plane with chromatic number at least $\log \log n$. In particular, it follows that the chromatic number of triangle-free string graphs cannot be bounded from above by a constant.

A {\it topological graph} is a graph drawn in the plane so that its vertices are represented by points and its edges are represented by (possibly crossing) curves connecting the corresponding point pairs. We also assume that no edge passes through any point representing a vertex other than its endpoints. For any integer $t\ge 2$,  we say that a topological graph is {\it $t$-quasi-planar} if it has no $t$ pairwise crossing edges. According to an old conjecture made independently by several people (see, e.g., Problem 6 in \cite{P91}), for any integer $t\ge 2$, there is a constant $c_t$ such that every $t$-quasi-planar topological graph on $n$ vertices has at most $c_tn$ edges. Any 2-quasi-planar graph is planar, therefore, it follows from Euler's polyhedral formula that the conjecture is true for $t=2$ with $c_2=3$. For $t=3$, extending earlier work of Agarwal, Aronov, Pach, Pollack, Sharir~\cite{AAPPS}, the conjecture was proved by Pach, Radoi\v ci\'c, and T\'oth \cite{PRT}. Subsequently, using the so-called ``discharging method" \cite{RT}, Ackerman \cite{A} also managed to prove the conjecture for $t = 4$. Theorem~\ref{thm1} immediately implies the following result, originally established in \cite{FP12b} by a more complicated argument.

\begin{corollary}\label{quasi} {\rm \cite{FP12b}}
Every $t$-quasi-planar topological graph on $n>2$ vertices has at most $n(\log n)^{c\log t}$ edges, for an appropriate constant $c$.
\end{corollary}

\noindent{\bf Proof} (using Theorem~\ref{thm1}): \hspace{1mm} Consider a $t$-quasi-planar graph $G$ with $n$ vertices and $m$ edges. Remove the endpoints of the edges, and take the intersection graph of the (open) edges as curves. We obtain a $K_t$-free string graph with $m$ vertices. By Theorem~\ref{thm1}, its vertices can be colored by at most $(\log m)^{C\log t}$ colors so that no two vertices of the same color are adjacent. Therefore, the curves corresponding to vertices of any given color class form a plane graph. Hence, the size of each color class is at most $3n$.  This implies that $m$, the total number of edges of $G$, satisfies the inequality $$m \leq 3n \cdot (\log m)^{C\log t} \leq 3n \cdot (2\log n)^{C\log t},$$ as required. {\hfill$\Box$\medskip}

A family of graphs is called {\em hereditary} if it is closed under induced subgraphs. The Erd\H{o}s-Hajnal conjecture \cite{EH89} states that for every hereditary family $F$ of graphs which is not the family of all graphs, there is a constant $c=c_F$ such that every graph in $F$ on $n$ vertices contains a clique or independent set of size $n^c$. A weaker estimate, with $e^{c\sqrt{\log n}}$ instead of $n^c$, was established by Erd\H{o}s and Hajnal.

The Erd\H{o}s-Hajnal conjecture is known to be true only for a few special classes of graphs; see the recent survey by Chudnovsky \cite{Ch} for partial results. Obviously, the families of intersection graphs of finitely many convex bodies, balls, curves, or other kinds of geometric objects in a given space, are hereditary. In many cases, it has been verified that these families satisfy the Erd\H os-Hajnal conjecture; see \cite{FP08b} for a survey. However, it is not known whether the Erd\H{o}s-Hajnal conjecture holds for string graphs. Our next theorem, represents the first progress on this problem. It is a lopsided statement: every string graph contains an independent set of size at least $n^c$ or a complete subgraph with at least $n^{c/\log \log n}$ vertices.

\begin{theorem}\label{thm2}
For every $\varepsilon>0$, there is a constant $c=c_{\varepsilon}>0$ such that the following holds. Every string graph on $n>2$ vertices contains a complete subgraph with at least $n^{c/\log \log n}$ vertices or an independent set of size $n^{1-\varepsilon}$. That is, every collection of $n>2$ curves in the plane contains a subcollection of at least $n^{c/\log \log n}$ pairwise intersecting curves or a subcollection of at least $n^{1-\varepsilon}$ pairwise disjoint curves.
\end{theorem}

\noindent {\bf Proof} (using Theorem \ref{thm1}): \hspace{1mm} Let $c=\frac{\varepsilon}{C}$. Applying Theorem~\ref{thm1} with $t=n^{c/\log \log n}$, we obtain that the chromatic number of any $K_t$-free string graph $G$ with $n$ vertices is at most $(\log n)^{C \log t}=t^{C \log\log n}=n^{\varepsilon}$. Thus, $G$ has an independent set of size at least $n/n^{\varepsilon}=n^{1-\varepsilon}$. {\hfill$\Box$\medskip}

The classical K\H ov\'ari-S\'os-Tur\'an theorem \cite{KST} states that any $K_{t,t}$-free graph with $n$ vertices has at most  $n^{2-1/t}+tn/2$ edges. Pach and Sharir~\cite{PS} conjectured that, for string graphs, this upper bound can be replaced by a bound linear in $n$. That is, every $K_{t,t}$-free string graph on $n$ vertices has at most $c_tn$ edges. They verified this conjecture up to a polylogarithmic factor in $n$. In \cite{FP10}, it was proved that the conjecture is true with $c_t \leq t^{c\log \log t}$. The authors further conjectured that the statement also holds with $c_t=ct \log t$, which would be best possible. We get close to this conjecture, proving the upper bound $c_t \leq t(\log t)^{O(1)}$.

\begin{theorem}\label{thm3}
There is a constant $c$ such that for any positive integers $t$ and $n$, every $K_{t,t}$-free string graph with $n$ vertices has at most $t(\log t)^c n$ edges.
\end{theorem}

The celebrated crossing lemma of Ajtai, Chv\'atal, Newborn, Szemer\'edi \cite{ACNS} and, independently, Leighton \cite{L84} states that in every drawing of a graph with $n$ vertices and $m \geq 4n$ edges, there are at least $\Omega(\frac{m^3}{n^2})$ pairs of crossing edges. This is easily seen to be equivalent to the existence of {\it one} edge that crosses  $\Omega(\frac{m^2}{n^2})$ other edges. Indeed, by the crossing lemma, the average number of edges a single edge crosses is $\Omega(\frac{m^2}{n^2})$. In the other direction, by repeatedly pulling out one edge at a time that crosses $\Omega(\frac{m^2}{n^2})$ of the remaining edges, a total of $\Omega(m)$ edges are pulled out that each cross $\Omega(\frac{m^2}{n^2})$ other edges. This gives $\Omega(m \cdot \frac{m^2}{n^2})$ pairs of crossing edges, and hence implies the crossing lemma.

Can the crossing lemma be strengthened to show that every graph drawn with $n$ vertices and $m \geq 4n$ edges contains two sets $E_1,E_2$ of edges, each of size $\Omega(\frac{m^2}{n^2})$, such that every edge in $E_1$ crosses every edge in $E_2$? In \cite{FPT}, the authors and Cs.~T\'oth proved that, although the answer is no, the statement is true up to a polylogarithmic factor. It is not hard to see that this result is an immediate consequence of Theorem~\ref{thm3}.

\begin{corollary}\label{crossinglemma} {\rm \cite{FPT}}
In every topological graph $G$ with $n$ vertices and $m \geq 4n$ edges, there are two disjoint sets edges, each of cardinality at least $\frac{m^2}{n^2(\log \frac{m}{n})^c}$, such that every edge in one set crosses all edges in the other. Here $c>0$ is a suitable absolute constant.
\end{corollary}

\noindent{\bf Proof} (using Theorem~\ref{thm3} and the crossing lemma): \hspace{1mm} It follows from the crossing lemma that the intersection graph of the edges of $G$ with their endpoints deleted is a string graph $G'$ with $m$ vertices and $\Omega(\frac{m^3}{n^2})$ edges. Thus, the average degree of the vertices in $G'$ is $\Omega(\frac{m^2}{n^2})$. Theorem \ref{thm3} then implies that $G'$ contains a complete bipartite graph $K_{t,t}$ as a subgraph, with $t=\Omega(\frac{m^2}{n^2(\log \frac{m}{n})^c})$. The two vertex classes of this bipartite graph correspond to the desired pair of crossing edge sets in $G$. {\hfill$\Box$\medskip}

At least one logarithmic factor is needed in Corollary~\ref{crossinglemma}. In \cite{FPT}, we constructed topological graphs with $n$ vertices and $m\ge 4n$ edges for which the largest pair of crossing sets has cardinality $O(\frac{m^2}{n^2\log \frac{m}{n}})$.

In the next two sections, we prove Theorems~\ref{thm1} and~\ref{thm3}, respectively. Apart from Matou\v sek's separator theorem, Lemma~\ref{lemsep}, our other main tool will be the following result established in \cite{FP12}. It shows that every dense string graph of $n$ vertices contains a complete bipartite subgraph such that each of its vertex classes is of size nearly $\Omega(n)$.

\begin{lemma}\label{lem2} {\rm \cite{FP12}}
There is a constant $b$ such that every string graph with $n$ vertices and $\epsilon n^2$ edges contains a complete bipartite graph with parts of order at least $\epsilon^{b} \frac{n}{\log n}$.
\end{lemma}

As mentioned in \cite{FP12}, the construction in \cite{F} shows that the dependence on $n$ is tight, giving a string graph on $n$ vertices with edge density $1-o(1)$ and whose largest balanced complete bipartite graph has $O(n/\log n)$ vertices.

\section{Proof of Theorem \ref{thm1}}

To prove Theorem \ref{thm1}, it suffices to establish the following lemma.

\begin{lemma}\label{lemind}
There is an absolute constant $C$ such that every $K_t$-free string graph on $n>2$ vertices contains an independent set of size at least $n\left(\log n\right)^{-C\log t}$.
\end{lemma}

\noindent{\bf Proof of Theorem \ref{thm1}} (using Lemma~\ref{lemind}): \hspace{1mm} The statement is trivial for $t \leq 2$, so we may assume $t > 2$. We obtain a proper vertex coloring of the $K_t$-free string graph $G$ on $n$ vertices by repeatedly pulling out maximum indepedent sets, and giving a new color to the elements of each such independent set. It follows from the lower bound on the independence number, given in Lemma \ref{lemind}, that after using at most $$\frac{n}{(n/2)\left(\log (n/2)\right)^{-C\log t}} \leq 2\left(\log n\right)^{C\log t}$$ different colors, at least half of the vertices of $G$ have been colored. Therefore, one can properly color all the vertices, using at most $$\sum_{i=0}^{\log n}2\left(\log (n/2^i)\right)^{C\log t} \leq \left(1+\log n\right) \cdot 2\left(\log n\right)^{C\log t} \leq  4\left(\log n\right)^{C\log t+1}$$  colors. This completes the proof of Theorem \ref{thm1}. Note that the constant $C$ in Theorem \ref{thm1} is a bit larger than the constant $C$ from Lemma \ref{lemind}. {\hfill$\Box$\medskip}

\noindent{\bf Proof of Lemma~\ref{lemind}:} \hspace{1mm}
Let $d \geq 1$ and $b$ be the constants from Lemmas \ref{lemsep} and \ref{lem2}, respectively. We may assume that $t > 2$, as otherwise the result is trivial. Let $C=\max(8d,6b+1)$. Let $I_t(n)$ denote the maximum $\alpha$ such that every $K_t$-free string graph $G=(V,E)$ on $n$ vertices contains an independent set of size $\alpha$. We will prove by induction on $n$ and $t$ that
\begin{equation}\label{eq1} I_t(n) \geq n\left(\log n\right)^{-C\log t},
\end{equation}
for every $n>2$.

In the base cases $3 \leq n \leq  2^{4d}$, we know that $G$ has an independent set of size $1 \geq n(\log n)^{-C\log t}$, where we used $C \geq 8d$ and $n,t \geq 3$. The proof splits into two cases, depending on whether the number $m$ of edges of $G$ is small or large. Let $\epsilon=\left(4d (\log n)^2\right)^{-2}$. As we may assume $n > 2^{4d}$, we have $\epsilon \geq (\log n)^{-6}$.

{\bf Case 1:} $m \leq \epsilon n^2$. In this case, by Lemma \ref{lemsep}, the string graph $G$ has a separator $S$ of size at most $d \sqrt{m}\log m \leq 2d\epsilon^{1/2} n \log n = n/\left(2\log n\right)$. Thus, there is a partition $V=S \cup V_1 \cup V_2$ with $|V_1|,|V_2| \leq 2n/3$ such that no vertex in $V_1$ is adjacent to any vertex in $V_2$. Let $n_i=|V_i|$ for $i=1,2$. The union of the largest independent set in $V_1$ and the largest independent set in $V_2$ is an independent set. Hence,
\begin{eqnarray*} I_t(n) & \geq & I_t(n_1)+I_t(n_2) \geq  n_1\left(\log n_1\right)^{-C\log t}+n_2\left(\log n_2\right)^{-C\log t} \geq  (n_1+n_2)\left(\log 2n/3\right)^{-C\log t}  \\ & \geq &
n\left(1-\frac{1}{2\log n}\right)\left(\log n\right)^{-C\log t} \left(1-\frac{\log 2/3}{\log n}\right)^{-C\log t} \\ & \geq & n\left(\log n\right)^{-C\log t} \left(1-\frac{1}{2\log n}\right)\left(1-\frac{\log 2/3}{\log n}\right)^{-1} \geq n\left(\log n\right)^{-C\log t},
\end{eqnarray*}
where the second inequality uses the induction hypothesis. This completes the proof in this case.

{\bf Case 2:} $m > \epsilon n^2$. In this case, by Lemma \ref{lem2}, the string graph $G$ has a complete bipartite graph with parts $A$ and $B$ of size $s \geq \epsilon^{b}n/\log n = n/(\log n)^{6b+1} \geq n(\log n)^{-C}$. As $G$ is $K_t$-free, and it contains all edges between $A$ and $B$, at least one of the subgraphs induced by $A$ or by $B$ is $K_{t/2}$-free. Thus,
\begin{eqnarray*} I_t(n) & \geq & I_{t/2}(s) \geq s(\log s)^{-C\log (t/2)} \geq s(\log n)^{-C\log (t/2)} = s(\log n)^{C-C\log t} \\ & \geq &  n(\log n)^{-C}(\log n)^{C-C\log t}=n(\log n)^{-C\log t},
\end{eqnarray*}
where the second inequality uses the induction hypothesis. This completes the proof.
{\hfill$\Box$\medskip}

\section{Proof of Theorem \ref{thm3}}

The proof relies on the following technical lemma from \cite{FP08}, whose proof is based on a simple divide and conquer approach. It shows that if every member of a hereditary family of graphs admits a small separator, then the number of edges of each graph in this family is at most linear in the number of vertices. Given a nonnegative function $f$ defined on the set of positive integers, we say that a family $F$ of graphs is {\it $f$-separable} if every graph in $F$ with $n$ vertices has a separator of size at most $f(n)$.

\begin{lemma}\label{fewedgessep} {\rm \cite{FP08}}
Let $n_0 \geq 1$ and $\phi(n)$ be a monotone decreasing nonnegative function for $n \geq n_0$ with $\phi(n_0) \leq \frac{1}{12}$. If $F$ is an $n\phi(n)$-separable hereditary family of graphs, then every graph in $F$ with $n \geq n_0$ vertices has fewer than $\frac{qn_0}{2}n$ edges, where
$$q=\prod_{i=0}^{\infty}\left(1+\phi((4/3)^i n_0\right).$$
\end{lemma}

Let $F$ be the hereditary family of graphs which consists of all $K_{t,t}$-free string graphs. It follows from Lemma \ref{lem2} that if $G$ is a $K_{t,t}$-free string graph with $n$ vertices, it must have fewer than $\epsilon n^2$ edges, where $t=\epsilon^{b}n/\log n$ and $b \geq 1$ is the constant from Lemma \ref{lem2}, so that we have $\epsilon=\left(\frac{t\log n}{n}\right)^{1/b}$. By Lemma \ref{lemsep}, $G$ has a separator of size at most $$d\sqrt{m}\log m \leq 2d\sqrt{m}\log n \leq 2d\epsilon^{1/2}n\log n =  n\phi(n),$$
where $d \geq 1$ is the constant that appears in Lemma \ref{lemsep} and $\phi(n)=2dt^{\frac{1}{2b}}(\log n)^{1+\frac{1}{2b}}n^{-\frac{1}{2b}}$. One can easily check by taking the derivative of $\phi$ that $\phi(n)$ is a monotone decreasing function for $n \geq e^{2b+1}$.

Let $n_0=xt(\log t)^a,$ where $x=(2^8db)^{16b}$ and $a=8b$. Then we have $n_0 \geq e^{2b+1}$, and hence $\phi(n)$ is a monotone decreasing function for $n \geq n_0$. We have \begin{eqnarray*}\phi(n_0) & = & 2d(t/n_0)^{\frac{1}{2b}}(\log n_0)^{1+\frac{1}{2b}}=2dx^{-\frac{1}{2b}}\left(\log t\right)^{-\frac{a}{2b}}\left( \log n_0\right)^{1+\frac{1}{2b}} \\ & \leq & 2^{-63}d^{-7}b^{-8}(\log t)^{-4}(\log n_0)^2 = 2^{-63}d^{-7}b^{-8}(\log t)^{-4}(\log x + \log t + 8b\log \log t)^2
\\ & \leq & 2^{-50}d^{-7}b^{-6}(\log t)^{-2}(\log x)^2 \leq 2^{-40}d^{-7}b^{-4}(\log t)^{-2}\left(\log (2^8db)\right)^2 \leq \frac{1}{12}.
\end{eqnarray*}

Note that, for $n \geq n_0$,
\begin{eqnarray*}\frac{\phi((4/3)n)}{\phi(n)} & = & \left(1+\frac{\log (4/3)}{\log n}\right)^{1+1/2b}(4/3)^{-\frac{1}{2b}} \leq \left(1+\frac{\log (4/3)}{\log n_0}\right)^{1+1/2b}(4/3)^{-\frac{1}{2b}} \\ & \leq & \left(1+\frac{2^{-8}}{b}\right)^{1+1/2b}(4/3)^{-\frac{1}{2b}} \leq \left(1+\frac{1}{100b}\right)(4/3)^{-\frac{1}{2b}} \leq 1-\frac{1}{12b}.\end{eqnarray*}

Hence, $$\prod_{i=0}^{\infty}(1+\phi\left((4/3)^i n_0\right) \leq \textrm{exp}\left(\sum_{i=0}^{\infty}\phi((4/3)^i n_0)\right) \leq \textrm{exp}\left(\sum_{i=0}^{\infty} \left(1-\frac{1}{12b}\right)^{-i}/12\right) =e^{b}.$$
Thus, by Lemma \ref{fewedgessep}, the number of edges of any graph in $F$ on $n$ vertices, i.e., any $K_{t,t}$-free string graph on $n$ vertices, is at most $\frac{e^{b}n_0}{2} n = O(t (\log t)^a n)$. This completes the proof of Theorem \ref{thm3}. \qed

\vspace{0.1cm} 

\noindent {\bf Final Remark:}
Recall that we conjectured that the logarithmic factor in Lemma \ref{lemsep} can be removed. 
\begin{conjecture}{\rm \cite{FP10}} \label{sepconj}
Every string graph with $m$ edges has a separator of order $O(\sqrt{m})$. 
\end{conjecture}

Conjecture \ref{sepconj} would imply improvements to the results in this paper. For example,  by modifying the proof of Theorem~\ref{thm3}, Conjecture \ref{sepconj} implies the following conjecture.
\begin{conjecture}{\rm \cite{FP12}}\label{extremalconj}
Every $K_{t,t}$-free string graph on $n$ vertices has at most $ct(\log t)n$ edges.
\end{conjecture}

This in turn implies the following conjectured improvement to Corollary \ref{crossinglemma}, which would be tight. 
\begin{conjecture}
In every topological graph with $n$ vertices and $m \geq 4n$ edges, there are two disjoint sets edges, each of cardinality $\Omega(\frac{m^2}{n^2(\log \frac{m}{n})})$, such that every edge in one set crosses all edges in the other. 
\end{conjecture}
 Conjecture \ref{sepconj} would also imply improved constants in Theorem \ref{thm1}, Corollary \ref{quasi}, and Theorem \ref{thm2}.

\vspace{0.1cm}

\noindent {\bf Acknowledgement:} We would like thank Jirka Matou\v sek for sharing and discussing his separator theorem for string graphs, which led to this note, as well as helpful comments on an early draft.

\end{document}